\documentclass[11pt]{amsart}
\usepackage{latexsym}
\usepackage{amsmath}
\usepackage{amssymb}

\newcommand{\eproof}{\mbox{\ }\hfill $\Box$ \par \vskip 10pt}

\newtheorem{Theorem}{Theorem}[section]
\newtheorem{lemma}[Theorem]{Lemma}
\newtheorem{prop}[Theorem]{Proposition}

\newtheorem{corol}[Theorem]{Corollary}

\numberwithin{equation}{section}

\def\cal{\mathcal}

\begin{document}

\title[Semiclassical resolvent estimates]{Semiclassical resolvent estimates for H\"older potentials}

\author[G. Vodev]{Georgi Vodev}

\address {Universit\'e de Nantes, Laboratoire de Math\'ematiques Jean Leray, 2 rue de la Houssini\`ere, BP 92208, 44322 Nantes Cedex 03, France}
\email{Georgi.Vodev@univ-nantes.fr}

\date{}

\begin{abstract} We first prove semiclassical resolvent estimates for the
Schr\"odinger operator in $\mathbb{R}^d$, $d\ge 3$, with real-valued potentials which are 
H\"older with respect to the radial variable. Then we extend these resolvent estimates 
to exterior domains in $\mathbb{R}^d$, $d\ge 2$, and real-valued potentials which are 
H\"older with respect to the space variable. As an application, we obtain the rate of the decay of the local energy
of the solutions to the wave equation with a refraction index which may be H\"older, Lipschitz or just $L^\infty$.

\quad

Key words: Schr\"odinger operator, resolvent estimates, H\"older potentials.
\end{abstract} 

\maketitle

\setcounter{section}{0}
\section{Introduction and statement of results}

In this paper we are going to 
study the resolvent of the Schr\"odinger operator
$$P(h)=-h^2\Delta+V(x)$$
where $0<h\le 1$ is a semiclassical parameter, $\Delta$ is the negative Laplacian in 
$\mathbb{R}^d$, $d\ge 2$, and $V\in L^\infty(\mathbb{R}^d)$ is a real-valued potential satisfying the condition
\begin{equation}\label{eq:1.1}
V(x)\le p(|x|)
\end{equation}
where $p(r)>0$, $r\ge 0$, is a decreasing function such that 
$p(r)\to 0$ as $r\to\infty$. More precisely, we are interested in bounding  
the quantity
$$g_s^\pm(h,\varepsilon):=\log\left\|(|x|+1)^{-s}(P(h)-E\pm i\varepsilon)^{-1}(|x|+1)^{-s}
\right\|_{L^2\to L^2}$$
from above by an explicit function of $h$, independent of $\varepsilon$, without imposing extra assumptions on the function $p$. 
Here $L^2:=L^2(\mathbb{R}^d)$, $0<\varepsilon<1$, $s>1/2$ is independent of $h$ and $E>0$ is a fixed energy level independent of $h$.
Instead, we impose some regularity on the potential with respect to the radial variable $r=|x|$. 
Note that througout this paper the space $C^1$ will denote the Lipschitz functions, that is, the ones with first derivatives
belonging to $L^\infty$ (and not necessairily continuous).

We will first extend Datchev's result \cite{kn:D} to a larger class of potentials. Recall that in \cite{kn:D}
the bound
\begin{equation}\label{eq:1.2}
g_s^\pm(h,\varepsilon)\le Ch^{-1}
\end{equation}
is proved when $d\ge 3$, with some constant $C>0$ independent of $h$ and $\varepsilon$, for potentials $V\in C^1(\overline{\mathbb{R}^+})$ 
with respect to the radial variable $r$ and satisfying (\ref{eq:1.1})
with $p(|x|)=C_1(|x|+1)^{-\delta}$ as well as the condition
\begin{equation}\label{eq:1.3}
\partial_rV(x)\le C_2(|x|+1)^{-\beta}
\end{equation}
where $C_1, C_2, \delta>0$ and $\beta>1$ are some constants. We will prove the following

\begin{Theorem} Let $d\ge 3$ and suppose that the potential $V$ satisfies the conditions (\ref{eq:1.1}) and (\ref{eq:1.3}). 
Then there exists a constant $C>0$ independent of $h$ and $\varepsilon$ but depending on $s$, $E$ and the function $p$, such 
that the bound (\ref{eq:1.2}) holds for all $0<h\le 1$.
\end{Theorem}

Note that the bound (\ref{eq:1.2}) was first proved for smooth potentials
in \cite{kn:B2}. 
A high-frequency analog of (\ref{eq:1.2}) on Riemannian manifolds was also proved in 
\cite{kn:B1} and \cite{kn:CV}. When $d=2$ the bound (\ref{eq:1.2}) is proved in \cite{kn:S2}
for potentials $V\in C^1(\mathbb{R}^2)$ satisfying (\ref{eq:1.1})
with $p(|x|)=C_1(|x|+1)^{-\delta}$ as well as the condition
\begin{equation}\label{eq:1.4}
|\nabla V(x)|\le C_2(|x|+1)^{-\beta}
\end{equation}
where $C_1, C_2, \delta>0$ and $\beta>1$ are some constants.

On the other hand, for compactly supported $L^\infty$ potentials without any regularity the following weaker bound
\begin{equation}\label{eq:1.5}
g_s^\pm(h,\varepsilon)\le Ch^{-4/3}\log(h^{-1})
\end{equation}
 was proved for $0<h\ll 1$ in \cite{kn:KV} and \cite{kn:S3} when $d\ge 2$. When $d\ge 3$ the bound (\ref{eq:1.5}) has been extended in
\cite{kn:V1} to potentials satisfying the condition
\begin{equation}\label{eq:1.6}
|V(x)|\le C_3(|x|+1)^{-\delta}
\end{equation}
where $C_3>0$ and $\delta>3$ are some constants. Note that (\ref{eq:1.5}) has been recently proved in \cite{kn:GS} for potentials satisfying
(\ref{eq:1.6}) with $\delta>2$. For potentials satisfying (\ref{eq:1.6})
with $1<\delta\le 3$ the much weaker bound
\begin{equation}\label{eq:1.7}
g_s^\pm(h,\varepsilon)\le 
 Ch^{-\frac{2\delta+5}{3(\delta-1)}}\left(\log(h^{-1})\right)^{\frac{1}{\delta-1}}.
\end{equation}
was proved in \cite{kn:V2}.

In the present paper we show that the bound (\ref{eq:1.5}) can be improved if some small regularity of
the potential is assumed. 
To be more precise, given $0<\alpha<1$ and $\beta>0$, 
we introduce the space $C_\beta^\alpha(\overline{\mathbb{R}^+})$ of all H\"older functions $a$ such that
 $$\sup_{r'\ge 0:\,0<|r-r'|\le 1}\frac{|a(r)-a(r')|}{|r-r'|^\alpha}\le C(r+1)^{-\beta},\quad\forall r\in \overline{\mathbb{R}^+},$$
 for some constant $C>0$. 
We now suppose that the function $V(r,w):=V(rw)$ satisfies the condition
\begin{equation}\label{eq:1.8}
V(\cdot,w)\in C_4^\alpha(\overline{\mathbb{R}^+}),\quad 0<\alpha<1,
\end{equation}
uniformly in $w\in \mathbb{S}^{d-1}$. 
We have the following

\begin{Theorem} Let $d\ge 3$ and suppose that the potential $V$ satisfies the conditions (\ref{eq:1.1}) and (\ref{eq:1.8}). 
Then there exists a constant $C>0$ independent of $h$ and $\varepsilon$ but depending on $s$, $E$ and the function $p$, such 
that the bound
\begin{equation}\label{eq:1.9}
g_s^\pm(h,\varepsilon)\le Ch^{-4/(\alpha+3)}\log(h^{-1})+C
\end{equation}
holds for all $0<h\le 1$.
\end{Theorem}

The proof of the above theorems is based on the global Carleman estimates proved in \cite{kn:V2} 
but with different phase and weight functions (see Theorem 4.1). In fact, in the case of H\"older or Lipschitz
potentials we need to construct better phase functions and hence get better Carleman estimates. 
Such functions are constructed in Section 2 modifying the construction in
\cite{kn:V2} in a suitable way. In order that the Carleman estimates (see (\ref{eq:4.1})
and (\ref{eq:4.6}) below) hold, 
the phase and weight functions must satisfy some inequalities (see (\ref{eq:2.5}), (\ref{eq:2.9})
and (\ref{eq:2.21}) below), so most of the proof of the above theorems consists of proving these
inequalities. Note also that the above theorems have been recently proved in \cite{kn:GS} by using similar Carleman
estimates but with a better choice of the phase function. Consequently, the bound (\ref{eq:1.9}) is proved in \cite{kn:GS} for a larger
class of $\alpha$-H\"older potentials. On the other hand, it is shown in \cite{kn:V4} that the logarithmic term in the right-hand side of
(\ref{eq:1.9}) can be removed for radial potentials. 

We next extend the above results to arbitrary obstacles and all dimensions $d\ge 2$. To do so, we need to replace the conditions 
(\ref{eq:1.3}) and (\ref{eq:1.8}) by stronger ones. To be more precise, we let $\Omega\subset \mathbb{R}^d$, $d\ge 2$, be a connected domain
with smooth boundary $\partial\Omega$ such that $\mathbb{R}^d\setminus\Omega$ is compact. 
Let $r_0>0$ be such that $\mathbb{R}^d\setminus\Omega\subset\{x\in\mathbb{R}^d:|x|\le r_0\}$. Given a real-valued potential
$V\in L^\infty(\Omega)$ satisfying (\ref{eq:1.1}) for $|x|\ge r_0$, we denote by $P(h)$ the Dirichlet self-adjoint
realisation of the operator $-h^2\Delta+V(x)$ on the Hilbert space $L^2(\Omega)$. We define the quantity $g_s^\pm$ in the same
way as above with
$L^2=L^2(\Omega)$. Given $0<\alpha\le 1$ and $\beta>0$, 
we introduce the space $C_\beta^\alpha(\overline{\Omega})$ of all H\"older functions $a$ such that
 $$\sup_{x'\in\overline{\Omega}:\,0<|x-x'|\le 1}\frac{|a(x)-a(x')|}{|x-x'|^\alpha}\le C(|x|+1)^{-\beta},\quad\forall x\in\overline{\Omega},$$
 for some constant $C>0$. Note that the case $\alpha=1$ corresponds to the Lipschitz functions. We suppose that
\begin{equation}\label{eq:1.10}
V\in C_\beta^\alpha(\overline{\Omega}),\quad 0<\alpha\le 1,\,\beta>1.
\end{equation}
We have the following

\begin{Theorem} Let $d\ge 2$ and suppose that the potential $V\in L^\infty(\Omega)$ satisfies (\ref{eq:1.1}) for $|x|\ge r_0$. If $V$ satisfies (\ref{eq:1.10})
with $\alpha=1$ and $\beta>1$, then the bound (\ref{eq:1.2}) holds for all $0<h\le 1$.
If $V$ satisfies (\ref{eq:1.10})
with $0<\alpha<1$ and $\beta=4$, then the bound (\ref{eq:1.9}) holds for all $0<h\le 1$.
\end{Theorem}

To prove this theorem we follow the same strategy as in \cite{kn:V3}, where the bound (\ref{eq:1.5}) is proved in all dimensions
$d\ge 2$ for potentials $V\in L^\infty(\Omega)$ satisfying (\ref{eq:1.6}). It consists of gluing up two different types of estimates - 
 one in a compact set coming from the local Carleman estimates proved in \cite{kn:LR}
 (see Theorem 3.1) with a global Carleman estimate outside a sufficently big compact (see Theorem 4.2). This is carried out in Section 4.
 
Theorem 1.3 together with Theorem 1.1 of \cite{kn:V3} allow us to get uniform bounds for the resolvent of 
the Dirichlet self-adjoint realisation, $G$, 
of the operator $-n(x)^{-1}\Delta$ in the Hilbert space $H=L^2(\Omega,n(x)dx)$, where
$n\in L^\infty(\Omega)$ is a real-valued function called
refraction index satisfying the conditions
\begin{equation}\label{eq:1.11}
n_1\le n(x)\le n_2\quad\mbox{in}\quad\Omega,
\end{equation}
with some constants $n_1,n_2>0$, and
\begin{equation}\label{eq:1.12}
|n(x)-1|\le C(|x|+1)^{-\delta}\quad\mbox{in}\quad\Omega,
\end{equation}
with some constants $C,\delta>0$. More precisely, we have the following

\begin{corol} Suppose that the function $n$ satisfies the conditions (\ref{eq:1.11}) and (\ref{eq:1.12}).
Then, given any $s>1/2$ and $\lambda_0>0$ there is a constant $C>0$ depending on $s$ and $\lambda_0$ such that
the estimate 
\begin{equation}\label{eq:1.13}
\|(|x|+1)^{-s}(G-\lambda^2\pm i\varepsilon)^{-1}(|x|+1)^{-s}\|_{H\to H}\le e^{C\psi(\lambda)}
\end{equation}
holds for all $\lambda\ge\lambda_0$ uniformly in $\varepsilon$, where $\psi(\lambda)=\lambda^{4/3}\log(\lambda+1)$ if $n\in L^\infty(\Omega)$ satisfies (\ref{eq:1.12})
with $\delta>3$, 
$\psi(\lambda)=\lambda^{4/(\alpha+3)}\log(\lambda+1)$ if $n\in C_4^\alpha(\overline\Omega)$ with $0<\alpha<1$,
$\psi(\lambda)=\lambda$ if $n\in C_\beta^1(\overline\Omega)$ with $\beta>1$.
\end{corol}

To get (\ref{eq:1.13}) we apply the theorems mentioned above with $h=\lambda_0/\lambda$, $V=\lambda_0^2(1-n)$, $E=\lambda_0^2$
and $\varepsilon$ replaced by $\varepsilon h^2n$.

Using Corollary 1.4 one can extend Shapiro's result \cite{kn:S1} on the local energy decay of the solutions of the following wave equation
\begin{equation}\label{eq:1.14}
\left\{
\begin{array}{lll}
(n(x)\partial_t^2-\Delta)u(t,x)=0\quad\mbox{in}\quad\mathbb{R}\times\Omega,\\
u(t,x)=0\quad\mbox{on}\quad \mathbb{R}\times\partial\Omega,\\
u(0,x)=f_1(x),\,\partial_tu(0,x)=f_2(x) \quad\mbox{in}\quad\Omega.
\end{array}
\right.
\end{equation}
Given any $r_0\gg 1$, denote $\Omega_{r_0}=\{x\in\Omega:|x|\le r_0\}$. 
We have the following

\begin{corol} Suppose that the function $n$ satisfies (\ref{eq:1.11}) and that $n=1$ outside some compact subset of $\Omega$. 
Then, the solution $u(t,x)$ to the equation (\ref{eq:1.14}) with compactly
supported initial data $(f_1,f_2)\in H_0^2(\Omega)\times H_0^1(\Omega)$
satisfies the estimate
\begin{equation}\label{eq:1.15}
\left\|\nabla u(t,\cdot)\right\|_{L^2(\Omega_{r_0})}+\left\|\partial_tu(t,\cdot)\right\|_{L^2(\Omega_{r_0})}\le 
C\omega(t)\left(\|f_1\|_{H^2(\Omega)}+\|f_2\|_{H^1(\Omega)}\right)
\end{equation}
for $t\gg 1$, where 
$$\omega(t)=\left(\frac{\log\log t}{\log t}\right)^{3/4}.$$
 Suppose in addition that  
$n\in C^\alpha(\overline\Omega)$ with $0<\alpha\le1$. Then the estimate (\ref{eq:1.15}) holds with
$$\omega(t)=\left(\frac{\log\log t}{\log t}\right)^{(\alpha+3)/4}$$
if $0<\alpha<1$, and with $\omega(t)=(\log t)^{-1}$ if $\alpha=1$. 
The estimate (\ref{eq:1.15}) remains valid when $\Omega=\mathbb{R}^d$.
\end{corol}

\noindent
{\bf Remark 1.} In view of the recent results in \cite{kn:V4}, when  $\Omega=\mathbb{R}^d$, $d\ge 3$ and the function
$n$ depends only on the radial variable $r$, the estimate (\ref{eq:1.15}) holds with $\omega(t)=(\log t)^{-3/4}$ if $n\in L^\infty$,
and with $\omega(t)=(\log t)^{-(\alpha+3)/4}$ if $n$ is $\alpha$ - H\"older in $r$.

Note that estimates similar to (\ref{eq:1.15}) were first proved by Burq \cite{kn:B1} in the case $n\equiv 1$. 
Note also that 
an analog of the above theorem is proved by Shapiro \cite{kn:S1} in the case $\Omega=\mathbb{R}^d$. Then an estimate similar to (\ref{eq:1.15}) is proved with 
$\omega(t)$ replaced by $(\log t)^{-3/4+\epsilon}$, $\epsilon>0$ being arbitrary. Moreover, if in addition the function $n$ is supposed Lipschitz,
then the decay rate is improved to $\omega(t)=(\log t)^{-1}$. The proof in \cite{kn:S1} is based on the resolvent estimates obtained in
\cite{kn:D}, \cite{kn:S2} and \cite{kn:S3}.

The assumption that $n=1$ outside some compact is only necessary to study the low-frequency behavior of the cut-off resolvent of $G$.
Indeed, under this assumption one can easily see that this behavior is exactly the same as in the case when $n\equiv 1$,
which in turn is well-known (e.g. see Appendix B.2 of \cite{kn:B1}). Therefore, in this case the low-frequency analysis can be carried out in precisely the same way as in \cite{kn:S1}. 
Most probably, the condition (\ref{eq:1.12}) with $\delta>2$ would be enough. The high-frequency analysis in our case is also very similar
to that one in \cite{kn:S1} with some slight modifications allowing to deduce from (\ref{eq:1.13}) the sharp decay rate $\omega(t)$
(instead of $(\log t)^{-3/4+\epsilon}$).

\section{Construction of the phase and weight functions} 

Let $\rho\in C_0^\infty([0,1])$, $\rho\ge 0$, be a real-valued function independent of $h$ such that $\int_0^\infty\rho(\sigma)d\sigma=1$.
If $V$ satisfies (\ref{eq:1.8}), we approximate it by the function
$$V_\theta(r,w)=\theta^{-1}\int_0^\infty\rho((r'-r)/\theta)V(r',w)dr'=\int_0^\infty\rho(\sigma)V(r+\theta\sigma,w)d\sigma$$
where $\theta=h^{2/(\alpha+3)}$. Indeed, we have
\begin{equation}\label{eq:2.1}
|V(r,w)-V_\theta(r,w)|\le\int_0^\infty\rho(\sigma)|V(r+\theta\sigma,w)-V(r,w)|d\sigma$$
$$\lesssim\theta^\alpha(r+1)^{-4}\int_0^\infty\sigma^\alpha\rho(\sigma)d\sigma\lesssim\theta^\alpha(r+1)^{-4}.
\end{equation}
This bound together with (\ref{eq:1.1}) imply
\begin{equation}\label{eq:2.2}
V_\theta(r,w)\le p(r)+{\cal O}((r+1)^{-4}).
\end{equation}
Clearly, $V_\theta$ is $C^1$ with respect to the variable $r$ and its first derivative $V'_\theta$ is given by
$$V'_\theta(r,w)=\theta^{-2}\int_0^\infty\rho'((r'-r)/\theta)V(r',w)dr'$$
$$=\theta^{-1}\int_0^\infty\rho'(\sigma)V(r+\theta\sigma,w)d\sigma=\theta^{-1}\int_0^\infty\rho'(\sigma)(V(r+\theta\sigma,w)-V(r,w))d\sigma$$
where we have used that $\int_0^\infty\rho'(\sigma)d\sigma=0$. Hence
\begin{equation}\label{eq:2.3}
|V'_\theta(r,w)|\lesssim\theta^{-1+\alpha}(r+1)^{-4}\int_0^\infty\sigma^\alpha|\rho'(\sigma)|d\sigma
\lesssim\theta^{-1+\alpha}(r+1)^{-4}.
\end{equation}
We now construct the weight function $\mu$ as follows:
$$\mu(r)=
\left\{
\begin{array}{lll}
 (r+1)^{2k}-(r+1)^{2k_0}&\mbox{for}& 0\le r\le a,\\
 (a+1)^{2k}-(a+1)^{2k_0}+(a+1)^{-2s+1}-(r+1)^{-2s+1}&\mbox{for}& r\ge a,
\end{array}
\right.
$$
where $a=a_0h^{-m}$ with $a_0\gg 1$ independent of $h$, $m=0$ if $V$ satisfies (\ref{eq:1.3})
and $m=2$ if $V$ satisfies (\ref{eq:1.8}). We choose $k=\frac{1}{4}\min\{1,\beta-1\}$, $k_0=0$ if $V$ satisfies (\ref{eq:1.3}),
and $k=1$, $k_0=1/2$ if $V$ satisfies (\ref{eq:1.8}). Furthermore, $s$ is independent of $h$ such that
\begin{equation}\label{eq:2.4}
\frac{1}{2}<s< 
\left\{
\begin{array}{lll}
\frac{1}{4}\min\{3,\beta+1\}& \mbox{if $V$ satisfies (\ref{eq:1.3})},\\
\frac{3}{4}& \mbox{if $V$ satisfies (\ref{eq:1.8})}.
\end{array}
\right.
\end{equation}
Clearly, the first derivative of $\mu$ is given by
$$
\mu'(r)=
\left\{
\begin{array}{lll}
 2k(r+1)^{2k-1}-2k_0(r+1)^{2k_0-1}&\mbox{for}& 0\le r<a,\\
 (2s-1)(r+1)^{-2s}&\mbox{for}& r>a.
\end{array}
\right.
$$
We have the following

\begin{lemma} For all $r>0$, $r\neq a$, we have the inequalities
\begin{equation}\label{eq:2.5}
2r^{-1}\mu(r)-\mu'(r)\ge 0,
\end{equation}
\begin{equation}\label{eq:2.6}
\frac{\mu(r)^j}{\mu'(r)}\lesssim a^{2kj}(r+1)^{2s},
\end{equation}
for every $j\ge 0$.
\end{lemma}

{\it Proof.} It is shown in Section 2 of \cite{kn:V2} that when $k_0=0$ the inequality (\ref{eq:2.5}) holds for all $0<k\le 1$.
Here we will prove it when $\nu:=2k-2k_0\ge 1$ and $0<k\le 1$. For $r<a$ we have
$$2\mu(r)-r\mu'(r)$$
$$=2(1-k)(r+1)^{2k}-2(1-k_0)(r+1)^{2k_0}+2k(r+1)^{2k-1}-2k_0(r+1)^{2k_0-1}$$
$$=2(r+1)^{2k_0-1}\left((1-k)(r+1)^{\nu+1}-(1-k_0)(r+1)+k(r+1)^{\nu}-k_0\right)$$
$$=2(r+1)^{2k_0-1}\left((1-k)r((r+1)^{\nu}-1)+(r+1)^{\nu}-\nu r/2-1\right)$$
$$\ge 2(r+1)^{2k_0-1}\left((r+1)^{\nu}-\nu r/2-1\right)\ge \nu r(r+1)^{2k_0-1}>0$$
where we have used the well-known inequality
$$(r+1)^{\nu}\ge\nu r+1$$
as long as $\nu\ge 1$. 
For $r>a$ the left-hand side of (\ref{eq:2.5}) is bounded
from below by
$$ 2r^{-1}((a+1)^{2k}-(a+1)^{2k_0}-s)>0$$
provided $a$ is taken large enough.  To prove (\ref{eq:2.6}) observe that for $r<a$ we have
$$\mu'(r)\ge 2(k-k_0)(r+1)^{2k-1}\ge 2(k-k_0)(r+1)^{-1}\ge 2(k-k_0)(r+1)^{-2s}$$
which clearly implies the bound (\ref{eq:2.6}) with $j=0$. This together with the fact that $\mu={\cal O}(a^{2k})$
implies the bound (\ref{eq:2.6}) with any $j>0$. 
\eproof

We will now construct a phase function
$\varphi\in C^1([0,+\infty))$ such that $\varphi(0)=0$ and $\varphi(r)>0$ for $r>0$. 
We define the first derivative of $\varphi$ by
$$\varphi'(r)=
\left\{
\begin{array}{lll}
 \tau(r+1)^{-k}- \tau(a+1)^{-k}&\mbox{for}& 0\le r\le a,\\
 0&\mbox{for}& r\ge a,
\end{array}
\right.
$$
where 
\begin{equation}\label{eq:2.7}
\tau=
\left\{
\begin{array}{lll}
\tau_0& \mbox{if $V$ satisfies (\ref{eq:1.3})},\\
\tau_0\theta^{2\alpha/3}h^{-1/3}& \mbox{if $V$ satisfies (\ref{eq:1.8})},
\end{array}
\right.
\end{equation}
with some parameter $\tau_0\gg 1$ independent of $h$ to be fixed later on. We choose now the parameter $a_0$ of the form
$a_0=\tau_0^\ell$, where $\ell>0$ is a constant such that $k\ell>2$ and $(\beta-2k-2s)\ell>2$. Note that the choice of the parameters
$k$ and $s$ guarantees that $\beta-2k-2s>0$. 

Clearly, the first derivative of $\varphi'$ satisfies
$$\varphi''(r)=
\left\{
\begin{array}{lll}
 -k\tau(r+1)^{-k-1}&\mbox{for}& 0\le r<a,\\
 0&\mbox{for}& r>a.
\end{array}
\right.
$$

\begin{lemma} For all $r\ge 0$ we have the bounds
\begin{equation}\label{eq:2.8}
h^{-1}\varphi(r)\lesssim 
\left\{
\begin{array}{lll}
 h^{-1}& \mbox{if $V$ satisfies (\ref{eq:1.3})},\\
 h^{-4/(\alpha+3)}\log(h^{-1})+1& \mbox{if $V$ satisfies (\ref{eq:1.8})},\\
\end{array}
\right.
\end{equation}
\end{lemma}

{\it Proof.} The lemma follows from the bounds
$$\max\varphi=\int_0^a\varphi'(r)dr\le \tau\int_0^a (r+1)^{-k}dr
\lesssim 
\left\{
\begin{array}{lll}
\tau a^{1-k}& \mbox{if}\quad k<1,\\
\tau\log a& \mbox{if}\quad k=1.\\
\end{array}
\right.
$$
\eproof

For $r>0$, $r\neq a$, set
$$A(r)=\left(\mu\varphi'^2\right)'(r),$$
$$B(r)=B_1(r)+B_2(r),$$
where
$$B_1(r)=(r+1)^{-\beta}\mu(r)+p(r)\mu'(r),$$
$$B_2(r)=\frac{\left(\mu(r)\varphi''(r)\right)^2}{h^{-1}\varphi'(r)\mu(r)+\mu'(r)},$$
with $\beta>1$, if $V$ satisfies (\ref{eq:1.3}), and 
$$B_1(r)=\theta^{-1+\alpha}(r+1)^{-\beta}\mu(r)+(p(r)+(r+1)^{-\beta})\mu'(r),$$
$$B_2(r)=\frac{\left(\mu(r)\left(h^{-1}\theta^\alpha(r+1)^{-\beta}+|\varphi''(r)|\right)\right)^2}{h^{-1}\varphi'(r)\mu(r)+\mu'(r)},$$
with $\beta=4$, if $V$ satisfies (\ref{eq:1.8}).
The following lemma will play a crucial role in the proof of the Carleman estimates (\ref{eq:4.1}) and (\ref{eq:4.6}) in the case
$d\ge 3$.

\begin{lemma} Given any constant $C>0$ there exists a positive constant $\tau_1=\tau_1(C,E)$ 
 such that for $\tau$ satisfying (\ref{eq:2.7}) with $\tau_0\ge\tau_1$ and for all $0<h\le 1$ we have the inequality
\begin{equation}\label{eq:2.9}
A(r)-CB(r)\ge -\frac{E}{2}\mu'(r)
\end{equation}
for all $r>0$, $r\neq a$. 
\end{lemma}

{\it Proof.} For $r<a$ we have
$$A(r)=-\left((r+1)^{2k_0}\varphi'^2\right)'+\tau^2\partial_r\left(1-(r+1)^k(a+1)^{-k}\right)^2$$ 
$$=-2(r+1)^{2k_0}\varphi'(r)\varphi''(r)-2k_0(r+1)^{2k_0-1}\varphi'(r)^2$$
$$-2k\tau^2(r+1)^{k-1}(a+1)^{-k}\left(1-(r+1)^k(a+1)^{-k}\right)$$
$$\ge 2\tau(k-k_0)(r+1)^{2k_0-k-1}\varphi'(r)-2k\tau^2(r+1)^{k-1}(a+1)^{-k}$$
$$\ge 2\tau(k-k_0)(r+1)^{2k_0-k-1}\varphi'(r)-{\cal O}\left(\tau^2 a^{-k}\right)\mu'(r)$$
$$\ge 2\tau(k-k_0)(r+1)^{2k_0-k-1}\varphi'(r)-{\cal O}\left(\tau_0^2 a_0^{-k}\right)\mu'(r)$$
$$\ge 2\tau(k-k_0)(r+1)^{2k_0-k-1}\varphi'(r)-{\cal O}\left(\tau_0^{-k\ell+2}\right)\mu'(r).$$
Hence, taking $\tau_0$ large enough, we can arrange the inequality
\begin{equation}\label{eq:2.10}
A(r)\ge 2\tau(k-k_0)(r+1)^{2k_0-k-1}\varphi'(r)-\frac{E}{4}\mu'(r)
\end{equation}
for all $r<a$.  
Observe now that if $0<r\le a/2$, then
\begin{equation}\label{eq:2.11}
\varphi'(r)\ge \gamma\tau(r+1)^{-k}
\end{equation}
with some constant $\gamma>0$. By (\ref{eq:2.10}) and (\ref{eq:2.11}) we conclude
\begin{equation}\label{eq:2.12}
A(r)\ge \widetilde\gamma\tau^2(r+1)^{-2(k-k_0)-1}-\frac{E}{4}\mu'(r)
\end{equation}
for all $r\le a/2$ with some constant $\widetilde\gamma>0$, and
\begin{equation}\label{eq:2.13}
A(r)\ge -\frac{E}{4}\mu'(r)\quad \mbox{for all}\quad r\neq a.
\end{equation}
We will now bound the function $B_1$ from above. Since the function $p$ is decreasing, tending to zero, there is $b>0$ such that
$$p(r)+(r+1)^{-\beta}\le \frac{E}{9C}\quad\mbox{for}\quad r\ge b.$$
Hence, for every $N>0$ there is a constant $C_N>0$ such that we have
\begin{equation}\label{eq:2.14}
(p(r)+(r+1)^{-\beta})\mu'(r)\le C_N(r+1)^{-N}+\frac{E}{9C}\mu'(r)\quad\mbox{for all}\quad r\neq a.
\end{equation}
Let $0<r<a$. Then $\mu(r)<(r+1)^{2k}$, and in view of (\ref{eq:2.14}) with $N$ big enough, we have
$$B_1(r)\le \widetilde C(r+1)^{2k-\beta}+\frac{E}{9C}\mu'(r),$$
if $V$ satisfies (\ref{eq:1.3}), and 
$$B_1(r)\le \widetilde C\theta^{-1+\alpha}(r+1)^{2k-\beta}+\frac{E}{9C}\mu'(r),$$
with $\beta=4$, if $V$ satisfies (\ref{eq:1.8}). Observe now that the choice of the parameters $k,k_0$ and $\theta$ guarantees that
$\beta-2k\ge 2(k-k_0)+1$ and $\theta^{-1+\alpha}=\theta^{4\alpha/3}h^{-2/3}$. Therefore, the above inequalities imply
\begin{equation}\label{eq:2.15}
B_1(r)\le {\cal O}\left(\tau_0^{-2}\right)\tau^2(r+1)^{-2(k-k_0)-1}+\frac{E}{9C}\mu'(r)\quad\mbox{for}\quad r\le a/2
\end{equation}
in both cases. Similarly, we get
\begin{equation}\label{eq:2.16}
B_1(r)\le {\cal O}\left(\tau^2a^{-\beta+1}\right)\mu'(r)+\frac{E}{9C}\mu'(r)\quad\mbox{for}\quad  a/2<r<a
\end{equation}
and
\begin{equation}\label{eq:2.17}
B_1(r)\le {\cal O}\left(\tau^2a^{-\beta+2k+2s}\right)\mu'(r)+\frac{E}{9C}\mu'(r)\quad\mbox{for}\quad  r>a.
\end{equation}
Since 
$$\tau^2 a^{-\beta+1}<\tau^2a^{-\beta+2k+2s}\le \tau_0^2a_0^{-\beta+2k+2s}=\tau_0^{-(\beta-2k-2s)\ell+2},$$
 we obtain from (\ref{eq:2.16}) and (\ref{eq:2.17}),
\begin{equation}\label{eq:2.18}
B_1(r)\le \frac{E}{8C}\mu'(r)\quad\mbox{for}\quad  r>a/2,\,r\neq a,
\end{equation}
provided $\tau_0$ is taken large enough.

We will now bound the function $B_2$ from above.  
We will first consider the case when 
$V$ satisfies (\ref{eq:1.8}). Let $0<r\le a/2$. In view of (\ref{eq:2.11}), we have
$$B_2(r)\lesssim \frac{\mu(r)\left(h^{-2}\theta^{2\alpha}(r+1)^{-2\beta}+\varphi''(r)^2\right)}{h^{-1}\varphi'(r)}$$ 
 $$\lesssim  h^{-1}\theta^{2\alpha}\frac{\mu(r)(r+1)^{-2\beta}}{\varphi'(r)}
 + h\frac{\mu(r)\varphi''(r)^2}{\varphi'(r)}$$ 
$$\lesssim \tau^{-1}\theta^{2\alpha}h^{-1}(r+1)^{3k-2\beta}+h\tau(r+1)^{k-2}$$
$$\lesssim \tau_0^{-3}\tau^2(r+1)^{-2(k-k_0)-1}+ \tau(r+1)^{k-2}$$
where we have used that $5k-2k_0<2\beta-1$. Since $3k-2k_0-1>0$, we have the inequality
$$(r+1)^{k-2}\le b^{3k-2k_0-1}(r+1)^{-2(k-k_0)-1}+b^{-k-1}(r+1)^{2k-1}$$
for every $b>1$. We take $b$ such that $b^{3k-2k_0-1}=b_0\tau$, where $b_0>0$ is a small parameter independent
of $\tau$ and $h$ to be fixed below. Then the above inequality takes the form
$$\tau(r+1)^{k-2}\lesssim b_0\tau^2(r+1)^{-2(k-k_0)-1}+\tau^{-\frac{2(1-k+k_0)}{3k-2k_0-1}}\mu'(r)$$
$$\lesssim b_0\tau^2(r+1)^{-2(k-k_0)-1}+\tau_0^{-1}\mu'(r).$$
Thus, taking $\tau_0$ big enough depending on $b_0$, $E$ and $C$, we get the bound
\begin{equation}\label{eq:2.19}
B_2(r)\le {\cal O}\left(\tau_0^{-1}+b_0\right)\tau^2(r+1)^{-2(k-k_0)-1}+\frac{E}{8C}\mu'(r)\quad\mbox{for}\quad 0<r\le a/2.
\end{equation}
When $V$ satisfies (\ref{eq:1.3}) we have $3k-2k_0-1\le 0$, and hence
$$\tau(r+1)^{k-2}\le\tau(r+1)^{-2(k-k_0)-1}\le\tau_0^{-1}\tau^2(r+1)^{-2(k-k_0)-1}.$$
Therefore, the inequality (\ref{eq:2.19}) still holds in this case.

 Let us now see that
\begin{equation}\label{eq:2.20}
B_2(r)\le \frac{E}{8C}\mu'(r)\quad\mbox{for}\quad r>a/2,\,r\neq a.
\end{equation}
Let $\frac{a}{2}<r<a$. Since in this case $\mu(r)/\mu'(r)={\cal O}(r)$, we get the bound
$$B_2(r)\lesssim \left(\frac{\mu(r)}{\mu'(r)}\right)^2\left(h^{-1}\theta^\alpha(r+1)^{-\beta}+|\varphi''(r)|\right)^2\mu'(r)$$
$$\lesssim \left(h^{-2}\theta^{2\alpha}(r+1)^{2-2\beta}+\tau^2(r+1)^{-2k}\right)\mu'(r)$$
$$\lesssim \left(h^{-2}a^{2-2\beta}+\tau^2a^{-2k}\right)\mu'(r)$$
$$\lesssim \left(h^{2m(\beta-1)-2}a_0^{2-2\beta}+h^{2m-2/3}\tau_0^2a_0^{-2k}\right)\mu'(r)$$
$$\lesssim \left(a_0^{2-2\beta}+\tau_0^2a_0^{-2k}\right)\mu'(r)\lesssim \left(\tau_0^{-2\ell(\beta-1)}+\tau_0^{-2k\ell+2}\right)\mu'(r)$$
which clearly implies (\ref{eq:2.20}) in this case, provided $\tau_0$ is taken big enough. 
Let $r>a$. Using (\ref{eq:2.6}) with $j=1$, we get
$$B_2(r)\lesssim \left(\frac{\mu(r)}{\mu'(r)}\right)^2\left(h^{-1}\theta^\alpha(r+1)^{-\beta}\right)^2\mu'(r)$$ 
$$\lesssim h^{-2}a^{4k}(r+1)^{4s-2\beta}\mu'(r)$$
$$\lesssim h^{-2}a^{4k+4s-2\beta}\mu'(r)$$
$$\lesssim h^{2m(\beta-2k-2s)-2}a_0^{4k+4s-2\beta}\mu'(r)$$
$$\lesssim a_0^{4k+4s-2\beta}\mu'(r)\lesssim\tau_0^{-2\ell(\beta-2k-2s)}\mu'(r)$$
which again implies (\ref{eq:2.20}), provided $\tau_0$ is taken big enough. Similarly, in the case when $V$ satisfies (\ref{eq:1.3})
one concludes that the inequality (\ref{eq:2.20}) holds for all $r>0$, $r\neq a$. 

It is easy to see that for $r\le a/2$ the estimate (\ref{eq:2.9}) follows from (\ref{eq:2.12}), (\ref{eq:2.15}) and (\ref{eq:2.19}) 
by taking $b_0$ and $\tau_0^{-1}$ small enough, while for 
$r\ge a/2$, $r\neq a$, it follows from (\ref{eq:2.13}), (\ref{eq:2.18}) and (\ref{eq:2.20}).
\eproof

\noindent
{\bf Remark 2.} It is easy to see from the proof that when $V$ satisfies (\ref{eq:1.8}) the inequality (\ref{eq:2.9}) holds
as long as $1/2\le k\le 1$, $k_0=k-1/2$. The choice $k=1$, however, provides the best resolvent bound in the semiclassical
regime, that is, for $0<h\le h_0$ with some constant $0<h_0\ll 1$. When $h_0<h\le 1$ the choice of $k$ does not really matter 
because in this case $g_s^\pm(h,\varepsilon)$ is upper bounded by a constant. For example, we may take $k=1/2$ and $k_0=0$.

The following lemmas will play a crucial role in the proof of the Carleman estimate (\ref{eq:4.6}) in the case $d=2$.

\begin{lemma} Given any constants $C,r_0>0$ there exists a positive constant $\tau_1=\tau_1(C,E,r_0)$ 
  such that for $\tau$ satisfying (\ref{eq:2.7}) with $\tau_0\ge\tau_1$ and for all $0<h\le h_0$, $0<h_0<1$ being a constant
  depending on $E$, $r_0$ and $\tau_0$, we have the inequality
\begin{equation}\label{eq:2.21}
A(r)-h^2r^{-3}\mu(r)-CB(r)\ge -\frac{2E}{3}\mu'(r)
\end{equation}
for all $r\ge r_0$, $r\neq a$. 
\end{lemma}

{\it Proof.} For $r_0\le r<a$ we have 
$$
h^2r^{-3}\mu(r)\lesssim h^2(r+1)^{-3}\mu(r)\lesssim h^2(r+1)^{-2}\mu'(r)\le \frac{E}{6}\mu'(r),$$
provided $h$ is taken small enough. 
For $r>a$, in view of (\ref{eq:2.6}) with $j=1$, we have
$$
h^2r^{-3}\mu(r)\lesssim h^2a^{2k}(r+1)^{2s-3}\mu'(r)\lesssim h^2a^{2k+2s-3}\mu'(r)$$
$$\lesssim h^{2-m(2k+2s-3)}a_0^{2k+2s-3}\mu'(r)\le \frac{E}{6}\mu'(r),$$
provided $h$ is taken small enough, depending on $a_0$. 
Clearly, (\ref{eq:2.21}) follows from these inequalities and (\ref{eq:2.9}).
\eproof

It is easy to see from the proof that when $V$ satisfies (\ref{eq:1.3}) the inequality (\ref{eq:2.21}) holds also for 
$h_0<h\le 1$. This is no longer true when $V$ satisfies (\ref{eq:1.8}) because in this case $2k+2s-3$ does not have the right sign. 
Therefore, to make (\ref{eq:2.21}) holds for $h$ not necessarily small, we need to make a new choice of the parameters $k$
and $k_0$ in order to change the sign of $2k+2s-3$ and for which Lemma 2.3 still holds. Thus, in view of Remark 2, in the
semiclassical regime ($0<h\le h_0$) we take $k=1$, $k_0=1/2$ and in the
classical regime ($h_0< h\le 1$) we take $k=1/2$, $k_0=0$. To cover the second case we need the following 

\begin{lemma} If $V$ satisfies (\ref{eq:1.8}) we take $k=1/2$
and $k_0=0$. Then, given any constants $C,r_0>0$ there exists a positive constant $\tau_1=\tau_1(C,E,r_0)$ 
  such that for $\tau$ satisfying (\ref{eq:2.7}) with $\tau_0\ge\tau_1$ the inequality (\ref{eq:2.21})
  holds for all $r\ge r_0$, $r\neq a$, and all $0<h\le 1$.  
\end{lemma}

{\it Proof.} For $r_0\le r\le a/2$ we have 
$$
h^2r^{-3}\mu(r)\lesssim (r+1)^{-3}\mu(r)\lesssim (r+1)^{-3+2k}\lesssim (r+1)^{-2(k-k_0)-1}.$$
For $a/2<r<a$ we have
$$h^2r^{-3}\mu(r)\lesssim (r+1)^{-2}\mu'(r)\lesssim a^{-2}\mu'(r)\lesssim a_0^{-2}\mu'(r)\le \frac{E}{6}\mu'(r),$$
provided $a_0$ is taken big enough. 
For $r>a$ we have
$$
h^2r^{-3}\mu(r)\lesssim a^{2k}(r+1)^{2s-3}\mu'(r)\lesssim a^{2k+2s-3}\mu'(r)\lesssim a_0^{2k+2s-3}\mu'(r)\le \frac{E}{6}\mu'(r),$$
provided $a_0$ is taken big enough. 
Then it is easy to see that (\ref{eq:2.21}) follows from these inequalities and Remark 2.
\eproof

\section{Carleman estimates for H\"older potentials on bounded domains} 

Throughout this section $X\subset \mathbb{R}^d$, $d\ge 2$, will be a bounded, connected domain with a smooth boundary $\partial X$. 
Introduce the operator
$$P(h)=-h^2\Delta+V(x)$$
where $0<h\le 1$ is a semiclassical parameter and $V\in L^\infty(X)$ is a real-valued potential. 
Let $U\subset X$, $U\neq\emptyset$, be an arbitrary open domain, independent of $h$, such that $\partial U\cap\partial X=\emptyset$ and let $z\in\mathbb{C}$,
$|z|\le C_0$, $C_0>0$ being a constant independent of $h$. We will also denote by $H_h^1$ the Sobolev space equipped with the semiclassical norm. 
Given any $0<\alpha\le 1$, denote by $C^\alpha(\overline X)$ the space of all functions $a$ such that
$$\|a\|_{C^\alpha}:=\sup_{x',x\in\overline{X}:\,0<|x-x'|\le 1}\frac{|a(x)-a(x')|}{|x-x'|^\alpha}<+\infty.$$
We have the following

\begin{Theorem} Let $V\in C^\alpha(\overline X)$ with $0<\alpha\le 1$. Then, there exists a positive constant $\gamma$ depending on $U$, $\|V\|_{C^\alpha}$ and $C_0$ but independent of $h$ such that 
for all $0<h\le 1$ we have the 
estimate
\begin{equation}\label{eq:3.1}
\|u\|_{H_h^1(X)}\le e^{\gamma h^{-4/(\alpha+3)}}\|(P(h)-z)u\|_{L^2(X)}+e^{\gamma h^{-4/(\alpha+3)}}\|u\|_{H_h^1(U)}
\end{equation}
for every $u\in H^2(X)$ such that $u|_{\partial X}=0$.
\end{Theorem}

It is proved in Section 2 of \cite{kn:V3} that for complex-valued potentials $V\in L^\infty(X)$ the estimate (\ref{eq:3.1}) holds with $\alpha=0$. The proof is based on 
 the local Carleman estimates proved in \cite{kn:LR}. We will follow the same strategy in the case of H\"older potentials as well. For such potentials we will get new local Carleman estimates by making use of the results of \cite{kn:LR}. To be more precise, we 
 let $W\subset X$ be a small open domain and let
$x$ be local coordinates in $W$. If $\Gamma:=\overline W\cap\partial X$ is not empty we choose $x=(x_1,x')$, $x_1>0$ being the normal coordinate in $W$ and $x'$ the tangential ones. Thus in these coordinates $\Gamma$ is given by $\{x_1=0\}$. Let $\phi,\phi_1\in C^\infty(\overline W)$
be real-valued functions such that  
${\rm supp}\,\phi\subset {\rm supp}\,\phi_1\subset \overline W$,  $\phi_1=1$ on ${\rm supp}\,\phi$.
When $V\in C^\alpha(\overline X)$ with $0<\alpha<1$ we approximate the function $\phi_1V$ by the smooth function
$$V_\theta(x)=\theta^{-1}\int_X \varrho((x'-x)/\theta)(\phi_1V)(x')dx'$$
where $\varrho\in C_0^\infty(|x|\le 1)$ is a real-valued function such that $\int_{\mathbb{R}^d}\varrho(x)dx=1$ and $0<\theta<1$ is a small parameter to be fixed later on.  
The fact that $V\in C^\alpha(\overline X)$ implies the bounds
\begin{equation}\label{eq:3.2}
|(\phi_1V)(x)-V_\theta(x)|\lesssim \theta^\alpha,
\end{equation}
\begin{equation}\label{eq:3.3}
|\partial_x^\beta V_\theta(x)|\lesssim \theta^{\alpha-1},
\end{equation}
for all multi-indices $\beta$ such that $|\beta|=1$. Set $\widetilde V=\theta^{1-\alpha}(V_\theta-z)$ if $V\in C^\alpha(\overline X)$ with $0<\alpha<1$, $\widetilde V=V-z$ if $V\in C^1(\overline X)$. In view of (\ref{eq:3.2}) and (\ref{eq:3.3}) we have $\partial_x^\beta\widetilde V(x)={\cal O}(1)$
uniformly in $\theta$, for all multi-indices $\beta$ such that $|\beta|\le 1$.

Let now $\psi\in C^\infty(\overline W)$ be a real-valued function independent of $h$ and $\theta$ 
such that 
\begin{equation}\label{eq:3.4}
\nabla\psi\neq 0\quad\mbox{in}\quad \overline W.
\end{equation}
If $\Gamma\neq\emptyset$ we also suppose that
\begin{equation}\label{eq:3.5}
\frac{\partial\psi}{\partial x_1}(0,x')>0\quad\mbox{for all}\quad x'.
\end{equation}
We set $\varphi=e^{\lambda\psi}$, where $\lambda>0$ is a big parameter to be fixed later on, independent of $h$ and $\theta$. 
Let $p(x,\xi)\in C^\infty(T^*W)$ be the principal symbol of the operator $-\Delta$ and let $0<\widetilde h\ll 1$ be a new semiclassical parameter. Then the principal symbol, $\widetilde p_\varphi$, of the operator 
$$e^{\varphi/\widetilde h}(-\widetilde h^2\Delta+\widetilde V)e^{-\varphi/\widetilde h}$$
is given by the formula
$$\widetilde p_\varphi(x,\xi)=p(x,\xi+i\nabla\varphi(x))+\widetilde V(x).$$
An easy computation shows that given any constant $C>0$ there is $\lambda=\lambda(C)$ such that the condition (\ref{eq:3.4})
for the function $\psi$ implies the following condition for the function $\varphi$:
\begin{equation}\label{eq:3.6}
\left\{{\rm Re}\,\widetilde p_\varphi, {\rm Im}\,\widetilde p_\varphi\right\}(x,\xi)\ge c_1\quad\mbox{for}\quad |\xi|\le C, 
\end{equation}
with some constant $c_1>0$ independent of $\theta$. On the other hand, if $C$ is taken large enough we can arrange the lower bound
\begin{equation}\label{eq:3.7}
\left|\widetilde p_\varphi(x,\xi)\right|\ge c_2|\xi|^2\quad\mbox{for}\quad |\xi|\ge C, 
\end{equation}
with some constant $c_2>0$ independent of $\theta$. 
If $\Gamma\neq\emptyset$ the condition (\ref{eq:3.5}) implies
\begin{equation}\label{eq:3.8}
\frac{\partial\varphi}{\partial x_1}(0,x')>0\quad\mbox{for all}\quad x'.
\end{equation}
Now we are in position to use Propositions 1 and 2 of \cite{kn:LR}, where the proof is based on the properties (\ref{eq:3.6}), (\ref{eq:3.7})
and (\ref{eq:3.8}).  
We have the following

\begin{prop} Let the function $u$ be as in Theorem 3.1. 
Then there exist constants $C_1, \widetilde h_0>0$ such that for all $0<\widetilde h\le\widetilde h_0$ we have the estimate
\begin{equation}\label{eq:3.9}
\int_X\left(|\phi u|^2+|\widetilde h\nabla(\phi u)|^2\right)e^{2\varphi/\widetilde h}dx\le 
C_1\widetilde h^{-1}\int_X|(-\widetilde h^2\Delta+\widetilde V)(\phi u)|^2e^{2\varphi/\widetilde h}dx.
\end{equation}
\end{prop}

We take $\widetilde h=h\theta^{(1-\alpha)/2}$ when $\alpha<1$ and we rewrite the inequality (\ref{eq:3.9}) as follows
$$\int_X\left(|\phi u|^2+\theta^{1-\alpha}|h\nabla(\phi u)|^2\right)e^{2\varphi/h\theta^{(1-\alpha)/2}}dx$$
$$\le C_1h^{-1}\theta^{3(1-\alpha)/2}\int_X|(-h^2\Delta+V_\theta-z)(\phi u)|^2e^{2\varphi/h\theta^{(1-\alpha)/2}}dx$$
 $$\le C_1h^{-1}\theta^{3(1-\alpha)/2}\int_X|(P(h)-z)(\phi u)|^2e^{2\varphi/h\theta^{(1-\alpha)/2}}dx$$
 $$+C_1h^{-1}\theta^{3(1-\alpha)/2}\sup|\phi_1V-V_\theta|^2\int_X|\phi u|^2e^{2\varphi/h\theta^{(1-\alpha)/2}}dx$$
 $$\le C_1h^{-1}\theta^{3(1-\alpha)/2}\int_X|(P(h)-z)(\phi u)|^2e^{2\varphi/h\theta^{(1-\alpha)/2}}dx$$
 $$+C_2h^{-1}\theta^{(3+\alpha)/2}\int_X|\phi u|^2e^{2\varphi/h\theta^{(1-\alpha)/2}}dx.$$
 We now take $\theta=h^{2/(\alpha+3)}\kappa^{2/(1-\alpha)}$, where $\kappa>0$ is a small parameter independent of $h$. 
Thus, taking $\kappa$ small enough we can absorb the last term in the right-hand side of the above inequality. When $\alpha=1$
we take $\widetilde h=h\kappa$. Thus we deduce from Proposition 3.2 the following

\begin{prop} Let the function $u$ be as in Theorem 3.1.
Then there exist constants $\widetilde C, \kappa_0>0$ such that for all $0<\kappa\le\kappa_0$ and all $0<h\le 1$ we have the estimate
\begin{equation}\label{eq:3.10}
\int_X\left(|\phi u|^2+|h\nabla(\phi u)|^2\right)e^{2\varphi/\kappa h^{4/(\alpha+3)}}dx$$
$$\le \widetilde C\kappa h^{-2(\alpha+1)/(\alpha+3)}\int_X|(P(h)-z)(\phi u)|^2e^{2\varphi/\kappa h^{4/(\alpha+3)}}dx.
\end{equation}
\end{prop}

Now Theorem 3.1 follows from Proposition 3.3 in precisely the same way as in Section 2 of \cite{kn:V3}, where the analysis is
carried out in the particular case $\alpha=0$. It is an easy observation that the general case requires no changes in the arguments,
 and therefore we omit the details. 
 
\section{Resolvent estimates}

The following global Carlemann estimate is similar to that one in Section 3 of \cite{kn:V2}
and can be proved in the same way. 
The proof will be carried out in Section 5. In what follows we set ${\cal D}_r=-ih\partial_r$.

\begin{Theorem} Let $d\ge 3$ and let the potential $V$ satisfy (\ref{eq:1.1}). Let also $V$ satisfy 
either (\ref{eq:1.3}) or (\ref{eq:1.8}) and let $s$ satisfy (\ref{eq:2.4}). Then, for all $0<h\le 1$, $0<\varepsilon\le 1$ and for all functions
$f\in H^2(\mathbb{R}^d)$ such that 
$$(|x|+1)^{s}(P(h)-E\pm i\varepsilon)f\in L^2(\mathbb{R}^d)$$
 we have the estimate 
 \begin{equation}\label{eq:4.1}
\|(|x|+1)^{-s}e^{\varphi/h}f\|_{L^2(\mathbb{R}^d)}+\|(|x|+1)^{-s}e^{\varphi/h}{\cal D}_rf\|_{L^2(\mathbb{R}^d)}$$
$$\le Ca^{2}h^{-1}\|(|x|+1)^{s}e^{\varphi/h}(P(h)-E\pm i\varepsilon)f\|_{L^2(\mathbb{R}^d)}$$ 
$$+C\tau a\left(\varepsilon/h\right)^{1/2}\|e^{\varphi/h}f\|_{L^2(\mathbb{R}^d)}
\end{equation}
with a constant $C>0$ independent of $h$, $\varepsilon$ and $f$.
\end{Theorem}

Theorems 1.1 and 1.2 can be obtained from Theorem 4.1 in the same way as in Section 4 of \cite{kn:V2}. 
We will sketch the proof for the sake of completeness. It follows from the estimate 
(\ref{eq:4.1}) and Lemma 2.2 that for $0<h\le 1$ and $s$ satisfying (\ref{eq:2.4}) we have the estimate 
\begin{equation}\label{eq:4.2}
\|(|x|+1)^{-s}f\|_{L^2}\le M\|(|x|+1)^{s}(P(h)-E\pm i\varepsilon)f\|_{L^2}
+M\varepsilon^{1/2}\|f\|_{L^2}
\end{equation}
where $M>0$ is given by 
$$\log M=\left\{
\begin{array}{lll}
 Ch^{-1}& \mbox{if $V$ satisfies (\ref{eq:1.3})},\\
 Ch^{-4/(\alpha+3)}\log(h^{-1})+C& \mbox{if $V$ satisfies (\ref{eq:1.8})},\\
\end{array}
\right.
$$
with a constant $C>0$ independent of $h$ and $\varepsilon$. On the other hand, since the operator $P(h)$ is symmetric, we have
$$
\varepsilon\|f\|^2_{L^2}=\pm{\rm Im}\,\langle (P(h)-E\pm i\varepsilon)f,f\rangle_{L^2}$$
$$\le (2M)^{-2}\|(|x|+1)^{-s}f\|^2_{L^2}+(2M)^2\|(|x|+1)^{s}(P(h)-E\pm i\varepsilon)f\|^2_{L^2}$$
which yields
\begin{equation}\label{eq:4.3}
M\varepsilon^{1/2}\|f\|_{L^2}\le \frac{1}{2}\|(|x|+1)^{-s}f\|_{L^2}+
2M^2\|(|x|+1)^{s}(P(h)-E\pm i\varepsilon)f\|_{L^2}.
\end{equation}
By (\ref{eq:4.2}) and (\ref{eq:4.3}) we get
\begin{equation}\label{eq:4.4}
\|(|x|+1)^{-s}f\|_{L^2}\le 4M^2\|(|x|+1)^{s}(P(h)-E\pm i\varepsilon)f\|_{L^2}.
\end{equation}
It follows from (\ref{eq:4.4}) that the resolvent estimate
\begin{equation}\label{eq:4.5}
\left\|(|x|+1)^{-s}(P(h)-E\pm i\varepsilon)^{-1}(|x|+1)^{-s}
\right\|_{L^2\to L^2}\le 4M^2
\end{equation}
holds for all $0<h\le 1$ and $s$ satisfying (\ref{eq:2.4}), and hence for all $s>1/2$ independent of $h$. Clearly,
(\ref{eq:4.5}) implies the desired bounds for $g_s^\pm$.

Given any $r_0>0$ we denote $Y_{r_0}:=\{x\in\mathbb{R}^d:|x|\ge r_0\}$ and we let $\eta_{r_0}\in C^\infty(\mathbb{R})$
be such that $\eta_{r_0}(r)=0$ for $r\le r_0/3$, $\eta_{r_0}(r)=1$ for $r\ge r_0/2$. We set $V_\eta(x):=\eta_{r_0}(|x|)V(x)$. 
To prove Theorem 1.3 we need the following

\begin{Theorem} Let $d\ge 3$ and let the potential $V$ satisfy (\ref{eq:1.1}) for $|x|\ge r_0$. Let also $V_\eta$ satisfy 
either (\ref{eq:1.3}) or (\ref{eq:1.8}) and let $s$ satisfy (\ref{eq:2.4}). Then, for all $0<h\le 1$, $0<\varepsilon\le 1$ and for all functions
$f\in H^2(Y_{r_0})$ such that $f=\partial_rf=0$ on $\partial Y_{r_0}$ and
$$(|x|+1)^{s}(P(h)-E\pm i\varepsilon)f\in L^2(Y_{r_0})$$
 we have the estimate 
 \begin{equation}\label{eq:4.6}
\|(|x|+1)^{-s}e^{\varphi/h}f\|_{L^2(Y_{r_0})}+\|(|x|+1)^{-s}e^{\varphi/h}{\cal D}_rf\|_{L^2(Y_{r_0})}$$
$$\le Ca^{2}h^{-1}\|(|x|+1)^{s}e^{\varphi/h}(P(h)-E\pm i\varepsilon)f\|_{L^2(Y_{r_0})}$$ $$
+C\tau a\left(\varepsilon/h\right)^{1/2}\|e^{\varphi/h}f\|_{L^2(Y_{r_0})}
\end{equation}
with a constant $C>0$ independent of $h$, $\varepsilon$ and $f$. 

Let $d=2$. If $V_\eta$ satisfies (\ref{eq:1.8}) and $k=1$, $k_0=1/2$, then
(\ref{eq:4.6}) holds for $0<h\le h_0$ with some constant $0<h_0\ll 1$ depending on $\tau_0$. 
If $V_\eta$ satisfies (\ref{eq:1.8}) and $k=1/2$, $k_0=0$, or $V_\eta$ satisfies (\ref{eq:1.3}), then
(\ref{eq:4.6}) holds for all $0<h\le 1$.
\end{Theorem}

The proof of Theorem 4.2 is similar to that one of Theorem 4.1 with some suitable modifications when $d=2$
and will be carried out in Section 5.

Theorem 1.3 can be derived from Theorems 3.1 and 4.2 in a way similar to the one developed in Section 5 of \cite{kn:V3}. 
Let $r_0>0$ be such that $Y_{r_0/3}\subset\Omega$. Fix 
$r_j$, $j=1,2,3,4$, such that $r_0<r_1<r_2<r_3<r_4$. Choose functions $\psi_1, \psi_2\in C^\infty(\mathbb{R}^d)$,
depending only on the radial variable $r$, such that
$\psi_1=1$ in $\mathbb{R}^d\setminus Y_{r_1}$, $\psi_1=0$ in $Y_{r_2}$, $\psi_2=1$ in $\mathbb{R}^d\setminus Y_{r_3}$, $\psi_2=0$ in $Y_{r_4}$. 
If $s$ satisfies (\ref{eq:2.4}), we choose a function $\chi_s\in C^\infty(\overline\Omega)$, $\chi_s>0$, such that $\chi_s(x)=|x|^{-s}$
on $Y_{r_0}$. 
Let $f\in H^2(\Omega)$ be such that $\chi_s^{-1}(P(h)-E\pm i\varepsilon)f\in L^2(\Omega)$
and $f|_{\partial\Omega}=0$. 
Set 
$${\cal Q}_0=\|\chi_s^{-1}(P(h)-E\pm i\varepsilon)f\|_{L^2(\Omega)},$$
$${\cal Q}_1=\|f\|_{L^2(Y_{r_1}\setminus Y_{r_2})}+\|{\cal D}_rf\|_{L^2(Y_{r_1}\setminus Y_{r_2})},$$
$${\cal Q}_2=\|f\|_{L^2(Y_{r_3}\setminus Y_{r_4})}+\|{\cal D}_rf\|_{L^2(Y_{r_3}\setminus Y_{r_4})},$$
and observe that
$$\|[P(h),\psi_j]f\|_{L^2}\lesssim {\cal Q}_j,\quad j=1,2.$$
We now apply Theorem 3.1 to the function $\psi_2f$ with $X=\Omega\setminus Y_{r_4}$ and $U\subset X$ such that $U\cap{\rm supp}\,\psi_2=
\emptyset$. Thus we obtain
\begin{equation}\label{eq:4.7}
\|f\|_{H_h^1(\Omega\setminus Y_{r_3})}\le \|\psi_2f\|_{H_h^1(\Omega\setminus Y_{r_4})}$$
$$\le e^{\gamma h^{-4/(\alpha+3)}}\|(P(h)-E\pm i\varepsilon)\psi_2f\|_{L^2(\Omega\setminus Y_{r_4})}$$
 $$\le e^{\gamma h^{-4/(\alpha+3)}}\|(P(h)-E\pm i\varepsilon)f\|_{L^2(\Omega\setminus Y_{r_4})}+e^{\gamma h^{-4/(\alpha+3)}}{\cal Q}_2
\end{equation}
with a constant $\gamma>0$ independent of $h$ and $\tau_0$. 
In particular, (\ref{eq:4.7}) implies
\begin{equation}\label{eq:4.8}
{\cal Q}_1\le e^{\gamma h^{-4/(\alpha+3)}}{\cal Q}_0+e^{\gamma h^{-4/(\alpha+3)}}{\cal Q}_2.
\end{equation}
On the other hand, it is clear that if $V$ satisfies (\ref{eq:1.10}) with $\alpha=1$ and $\beta>1$ (resp. $0<\alpha<1$ and $\beta=4$),
then $V_\eta$ satisfies (\ref{eq:1.3}) (resp. (\ref{eq:1.8})). Therefore, we can apply Theorem 4.2 to the function $(1-\psi_1)f$ to obtain 
\begin{equation}\label{eq:4.9}
\|(|x|+1)^{-s}e^{\varphi/h}f\|_{L^2(Y_{r_2})}+\|(|x|+1)^{-s}e^{\varphi/h}{\cal D}_rf\|_{L^2(Y_{r_2})}$$
 $$\le\|(|x|+1)^{-s}e^{\varphi/h}(1-\psi_1)f\|_{L^2(Y_{r_1})}+\|(|x|+1)^{-s}e^{\varphi/h}{\cal D}_r(1-\psi_1)f\|_{L^2(Y_{r_1})}$$
$$\le Ca^{2}h^{-1}\|(|x|+1)^{s}e^{\varphi/h}(P(h)-E\pm i\varepsilon)(1-\psi_1)f\|_{L^2(Y_{r_1})}$$ 
$$+C\tau a(\varepsilon/h)^{1/2}\|e^{\varphi/h}f\|_{L^2(Y_{r_1})}$$
$$\le Ca^{2}h^{-1}\|(|x|+1)^{s}e^{\varphi/h}(P(h)-E\pm i\varepsilon)f\|_{L^2(Y_{r_1})}
+Ca^{2}h^{-1}e^{\varphi(r_2)/h}{\cal Q}_1$$ 
$$+C\tau a(\varepsilon/h)^{1/2}\|e^{\varphi/h}f\|_{L^2(Y_{r_1})}
\end{equation}
for all $0<h\le 1$. 
In particular, (\ref{eq:4.9}) implies
\begin{equation}\label{eq:4.10}
e^{\varphi(r_3)/h}{\cal Q}_2\le Ca^{2}h^{-1}e^{\max\varphi/h}{\cal Q}_0+C\tau a(\varepsilon/h)^{1/2}
e^{\max\varphi/h}\|f\|_{L^2(\Omega)}$$
$$+Ca^{2}h^{-1}e^{\varphi(r_2)/h}{\cal Q}_1.
\end{equation}
We have 
$$\varphi(r_3)-\varphi(r_2)=\tau\int_{r_2}^{r_3}\left((r+1)^{-k}-(a+1)^{-k}\right)dr\ge c\tau$$
with some constant $c>0$.  We deduce from (\ref{eq:4.10})
\begin{equation}\label{eq:4.11}
{\cal Q}_2\le \exp\left(\widetilde\beta h^{-4/(\alpha+3)}+\max\varphi/h\right){\cal Q}_0$$
$$+\varepsilon^{1/2}
\exp\left(\widetilde\beta h^{-4/(\alpha+3)}+\max\varphi/h\right)\|f\|_{L^2(\Omega)}$$
$$+\tau_0^{2\ell}\exp\left((\beta-c\tau_0)h^{-4/(\alpha+3)}\right){\cal Q}_1
\end{equation}
with a constant $\widetilde\beta>0$ independent of $h$ and a constant $\beta>0$ independent of $h$ and $\tau_0$. Combining (\ref{eq:4.8}) and (\ref{eq:4.11}) we get
\begin{equation}\label{eq:4.12}
{\cal Q}_2\le \exp\left((\widetilde\beta+\gamma)h^{-4/(\alpha+3)}+\max\varphi/h\right){\cal Q}_0$$
$$+\varepsilon^{1/2}
\exp\left(\widetilde\beta h^{-4/(\alpha+3)}+\max\varphi/h\right)\|f\|_{L^2(\Omega)}$$
$$+\tau_0^{2\ell}\exp\left((\beta+\gamma-c\tau_0)h^{-4/(\alpha+3)}\right){\cal Q}_2.
\end{equation}
Taking $\tau_0$ big enough, independent of $h$, we can arrange that
$$\tau_0^{2\ell}\exp\left((\beta+\gamma-c\tau_0)h^{-4/(\alpha+3)}\right)\le \tau_0^{2\ell}\exp\left(-c\tau_0h^{-4/(\alpha+3)}/2\right)
\le \tau_0^{2\ell}\exp\left(-c\tau_0/2\right)\le 1/2$$
for all $0<h\le 1$. 
Thus we can absorb the last term in the right-hand side of (\ref{eq:4.12}) 
to conclude that 
\begin{equation}\label{eq:4.13}
{\cal Q}_1+{\cal Q}_2\le \exp\left(\beta_1 h^{-4/(\alpha+3)}+\max\varphi/h\right){\cal Q}_0$$
$$+\varepsilon^{1/2}
\exp\left(\beta_1 h^{-4/(\alpha+3)}+\max\varphi/h\right)\|f\|_{L^2(\Omega)}
\end{equation}
with a constant $\beta_1>0$ independent of $h$. By (\ref{eq:4.7}), (\ref{eq:4.9}) and (\ref{eq:4.13}) we obtain
\begin{equation}\label{eq:4.14}
\|\chi_sf\|_{L^2(\Omega)}\le N{\cal Q}_0+\varepsilon^{1/2}
N\|f\|_{L^2(\Omega)}
\end{equation}
where
$$N=\exp\left(\beta_2 h^{-4/(\alpha+3)}+\max\varphi/h\right)$$
with a constant $\beta_2>0$ independent of $h$. In the same way as above, using the fact that the operator $P(h)$ is symmetric, we get from (\ref{eq:4.14})
 that the resolvent estimate
\begin{equation}\label{eq:4.15}
\left\|\chi_s(P(h)-E\pm i\varepsilon)^{-1}\chi_s
\right\|_{L^2(\Omega)\to L^2(\Omega)}\le 4N^2
\end{equation}
holds for all $0<h\le 1$, $0<\varepsilon\le 1$ and $s$ satisfying (\ref{eq:2.4}), which together with Lemma 2.2 clearly imply the desired bound.

\section{Proof of Theorems 4.1 and 4.2}
 
The main point is to work with the polar coordinates $(r,w)\in\mathbb{R}^+\times\mathbb{S}^{d-1}$, $r=|x|$, $w=x/|x|$ and to use that $L^2(\mathbb{R}^d)=L^2(\mathbb{R}^+\times\mathbb{S}^{d-1}, r^{d-1}drdw)$. 
In what follows in this section we denote by $\|\cdot\|$ and $\langle\cdot,\cdot\rangle$
the norm and the scalar product in $L^2(\mathbb{S}^{d-1})$. We will make use of the identity
\begin{equation}\label{eq:5.1}
 r^{(d-1)/2}\Delta  r^{-(d-1)/2}=\partial_r^2+\frac{\widetilde\Delta_w}{r^2}
\end{equation}
where $\widetilde\Delta_w=\Delta_w-\frac{1}{4}(d-1)(d-3)$ and $\Delta_w$ denotes the negative Laplace-Beltrami operator
on $\mathbb{S}^{d-1}$. Set $u=r^{(d-1)/2}e^{\varphi/h}f$ and
$${\cal P}^\pm(h)=r^{(d-1)/2}(P(h)-E\pm i\varepsilon)r^{-(d-1)/2},$$
$${\cal P}^\pm_\varphi(h)=e^{\varphi/h}{\cal P}^\pm(h)e^{-\varphi/h}.$$
Using (\ref{eq:5.1}) we can write the operator ${\cal P}^\pm(h)$ in the coordinates $(r,w)$ as follows
$${\cal P}^\pm(h)={\cal D}_r^2+\frac{\Lambda_w}{r^2}-E\pm i\varepsilon +V$$
where we have put ${\cal D}_r=-ih\partial_r$ and $\Lambda_w=-h^2\widetilde\Delta_w$. Since the function $\varphi$
depends only on the variable $r$, we get
$${\cal P}^\pm_\varphi(h)={\cal D}_r^2+\frac{\Lambda_w}{r^2}-E\pm i\varepsilon -\varphi'^2+h\varphi''+
2i\varphi'{\cal D}_r+V.$$ 
We write $V=V_L+ V_S$ with 
$V_L:=V_\theta$ and $V_S:=V-V_\theta$ if $V$ satisfies (\ref{eq:1.8}), and $V_L:=V$ and $V_S:=0$ if $V$ satisfies (\ref{eq:1.3}). 
For $r>0$, $r\neq a$, introduce the function
$$F(r)=-\langle (r^{-2}\Lambda_w-E-\varphi'(r)^2+ V_L(r,\cdot))u(r,\cdot),u(r,\cdot)\rangle+\|{\cal D}_ru(r,\cdot)\|^2$$
where $ V_L(r,w):= V_L(rw)$. 
Then its first derivative is given by
$$F'(r)=\frac{2}{r}\langle r^{-2}\Lambda_wu(r,\cdot),u(r,\cdot)\rangle
+((\varphi')^2- V_L)'\|u(r,\cdot)\|^2$$
$$-2h^{-1}{\rm Im}\,\langle {\cal P}^\pm_\varphi(h)u(r,\cdot),{\cal D}_ru(r,\cdot)\rangle$$
$$\pm 2\varepsilon h^{-1}{\rm Re}\,\langle u(r,\cdot),{\cal D}_ru(r,\cdot)\rangle+4h^{-1}\varphi'\|{\cal D}_ru(r,\cdot)\|^2$$ 
$$+2h^{-1}{\rm Im}\,\langle (V_S+h\varphi'')u(r,\cdot),{\cal D}_ru(r,\cdot)\rangle.$$
Thus we obtain the identity
$$(\mu F)'=\mu'F+\mu F'$$
$$=(2r^{-1}\mu-\mu')\langle r^{-2}\Lambda_wu(r,\cdot),u(r,\cdot)\rangle$$ 
$$+(E\mu'+(\mu(\varphi')^2-\mu V_L)')\|u(r,\cdot)\|^2$$
$$-2h^{-1}\mu{\rm Im}\,\langle {\cal P}^\pm_\varphi(h)u(r,\cdot),{\cal D}_ru(r,\cdot)\rangle$$
$$\pm 2\varepsilon h^{-1}\mu{\rm Re}\,\langle u(r,\cdot),{\cal D}_ru(r,\cdot)\rangle+(\mu'+4h^{-1}\varphi'\mu)\|{\cal D}_ru(r,\cdot)\|^2$$ 
$$+2h^{-1}\mu{\rm Im}\,\langle (V_S+h\varphi'')u(r,\cdot),{\cal D}_ru(r,\cdot)\rangle.$$
Using that $\Lambda_w\ge 0$ as long as $d\ge 3$ together with (\ref{eq:2.5}) we get the inequality
$$\mu'F+\mu F'\ge (E\mu'+(\mu(\varphi')^2-\mu V_L)')\|u(r,\cdot)\|^2$$
$$+(\mu'+4h^{-1}\varphi'\mu)\|{\cal D}_ru(r,\cdot)\|^2$$
$$-\frac{3h^{-2}\mu^2}{\mu'}\|{\cal P}^\pm_\varphi(h)u(r,\cdot)\|^2-\frac{\mu'}{3}\|{\cal D}_ru(r,\cdot)\|^2$$
$$-\varepsilon h^{-1}\mu\left(\|u(r,\cdot)\|^2+\|{\cal D}_ru(r,\cdot)\|^2\right)$$
$$-3h^{-2}\mu^2(\mu'+4h^{-1}\varphi'\mu)^{-1}\|(V_S+h\varphi'')u(r,\cdot)\|^2$$ 
$$-\frac{1}{3}(\mu'+4h^{-1}\varphi'\mu)\|{\cal D}_ru(r,\cdot)\|^2$$
$$\ge (E\mu'+(\mu(\varphi')^2)'-T_L\mu-Z_L\mu')\|u(r,\cdot)\|^2$$
$$+\frac{\mu'}{3}\|{\cal D}_ru(r,\cdot)\|^2-\frac{3h^{-2}\mu^2}{\mu'}\|{\cal P}^\pm_\varphi(h)u(r,\cdot)\|^2$$
$$-\varepsilon h^{-1}\mu\left(\|u(r,\cdot)\|^2+\|{\cal D}_ru(r,\cdot)\|^2\right)$$
$$-3h^{-2}\mu^2(\mu'+4h^{-1}\varphi'\mu)^{-1}(Q_S+h|\varphi''|)^2\|u(r,\cdot)\|^2$$ 
where 
$$T_L={\cal O}\left((r+1)^{-\beta}\right),\quad Z_L=p(r),\quad Q_S=0,$$
 if $V$ satisfies (\ref{eq:1.3}), 
$$T_L={\cal O}\left(\theta^{-1+\alpha}(r+1)^{-4}\right),\quad Z_L=p(r)+{\cal O}\left((r+1)^{-4}\right),\quad Q_S={\cal O}\left(\theta^\alpha(r+1)^{-4}\right),$$
if $V$ satisfies (\ref{eq:1.8}), and we have used the bounds (\ref{eq:2.1}),(\ref{eq:2.2}) and (\ref{eq:2.3}) in the second case.
Hence we can rewrite the above inequality in the form
 $$\mu'F+\mu F'\ge\left(E\mu'+A(r)-CB(r)\right)\|u(r,\cdot)\|^2+\frac{\mu'}{3}\|{\cal D}_ru(r,\cdot)\|^2$$
$$-\frac{3h^{-2}\mu^2}{\mu'}\|{\cal P}^\pm_\varphi(h)u(r,\cdot)\|^2
-\varepsilon h^{-1}\mu\left(\|u(r,\cdot)\|^2+\|{\cal D}_ru(r,\cdot)\|^2\right)$$
with a suitable constant $C>0$. 
 Now we use Lemma 2.3 to conclude that
 \begin{equation}\label{eq:5.2}
\mu'F+\mu F'\ge \frac{E}{2}\mu'\|u(r,\cdot)\|^2+\frac{\mu'}{3}\|{\cal D}_ru(r,\cdot)\|^2-\frac{3h^{-2}\mu^2}{\mu'}\|{\cal P}^\pm_\varphi(h)u(r,\cdot)\|^2$$ $$
-\varepsilon h^{-1}\mu\left(\|u(r,\cdot)\|^2+\|{\cal D}_ru(r,\cdot)\|^2\right).
\end{equation}
We integrate this inequality with respect to $r$ and use that $\mu(0)=0$. We have
$$\int_0^\infty(\mu'F+\mu F')dr=0.$$
Thus we obtain the estimate
\begin{equation}\label{eq:5.3}
\frac{E}{2}\int_0^\infty\mu'\|u(r,\cdot)\|^2dr+\int_0^\infty\frac{\mu'}{3}\|{\cal D}_ru(r,\cdot)\|^2dr\le 3h^{-2}\int_0^\infty\frac{\mu^2}{\mu'}
\|{\cal P}^\pm_\varphi(h)u(r,\cdot)\|^2dr$$ $$
+\varepsilon h^{-1}\int_0^\infty\mu\left(\|u(r,\cdot)\|^2+\|{\cal D}_ru(r,\cdot)\|^2\right)dr.
\end{equation}
Using that $\mu={\cal O}(a^{2})$ together with (\ref{eq:2.6}) we get from (\ref{eq:5.3})
\begin{equation}\label{eq:5.4}
\int_0^\infty(r+1)^{-2s}\left(\|u(r,\cdot)\|^2+\|{\cal D}_ru(r,\cdot)\|^2\right)dr$$
$$\le Ca^{4}h^{-2}\int_0^\infty(r+1)^{2s}\|{\cal P}^\pm_\varphi(h)u(r,\cdot)\|^2dr$$ $$
+C\varepsilon h^{-1}a^{2}\int_0^\infty\left(\|u(r,\cdot)\|^2+\|{\cal D}_ru(r,\cdot)\|^2\right)dr
\end{equation}
with some constant $C>0$ independent of $h$ and $\varepsilon$. On the other hand, we have the identity
$${\rm Re}\,\int_0^\infty\langle 2i\varphi'{\cal D}_ru(r,\cdot),u(r,\cdot)\rangle dr=\int_0^\infty h\varphi''\|u(r,\cdot)\|^2dr$$
and hence
$${\rm Re}\,\int_0^\infty\langle {\cal P}^\pm_\varphi(h)u(r,\cdot),u(r,\cdot)\rangle dr
=\int_0^\infty\|{\cal D}_ru(r,\cdot)\|^2dr
+\int_0^\infty \langle r^{-2}\Lambda_wu(r,\cdot),u(r,\cdot)\rangle dr$$
$$-\int_0^\infty(E+\varphi'^2)\|u(r,\cdot)\|^2dr
+\int_0^\infty\langle Vu(r,\cdot),u(r,\cdot)\rangle dr$$
$$\ge \int_0^\infty\|{\cal D}_ru(r,\cdot)\|^2dr-{\cal O}(\tau^2)\int_0^\infty\|u(r,\cdot)\|^2dr.$$
This implies
\begin{equation}\label{eq:5.5}
\varepsilon h^{-1}a^{2}\int_0^\infty\|{\cal D}_ru(r,\cdot)\|^2dr\le {\cal O}(\tau^2)\varepsilon h^{-1}a^{2}\int_0^\infty\|u(r,\cdot)\|^2dr$$
$$+\gamma\int_0^\infty(r+1)^{-2s}\|u(r,\cdot)\|^2dr
+\gamma^{-1}h^{-2}a^{4}\int_0^\infty(r+1)^{2s}\|{\cal P}^\pm_\varphi(h)u(r,\cdot)\|^2dr
\end{equation}
for every $\gamma>0$. Taking $\gamma$ small enough, independent of $h$, $\tau$ and $a$, and combining
the estimates (\ref{eq:5.4}) and (\ref{eq:5.5}), we get
\begin{equation}\label{eq:5.6}
\int_0^\infty(r+1)^{-2s}\left(\|u(r,\cdot)\|^2+\|{\cal D}_ru(r,\cdot)\|^2\right)dr$$
$$\le Ca^{4}h^{-2}\int_0^\infty(r+1)^{2s}\|{\cal P}^\pm_\varphi(h)u(r,\cdot)\|^2dr$$ $$
+C\varepsilon h^{-1}a^{2}\tau^2\int_0^\infty\|u(r,\cdot)\|^2dr
\end{equation}
with a new constant $C>0$ independent of $h$ and $\varepsilon$. Clearly, the estimate 
(\ref{eq:5.6}) implies (\ref{eq:4.1}).

The proof of Theorem 4.2 in the case when $d\ge 3$ goes very much like the proof of Theorem 4.1 above. The only difference in this case
is that we have to integrate the function $F(r)$ from $r_0$ to $\infty$ and use that $F(r_0)=0$ by assumption. 
Thus, by Lemma 2.3 we conclude that the inequality (\ref{eq:5.2}) holds for all $r\ge r_0$.

In the case $d=2$ the operator $\Lambda_w$ is no longer non-negative. Instead, we will use that so is the operator
$-\Delta_w$. Thus, it is easy to see that the above inequalities still hold with $V_L$ replaced by $V_L-h^2(2r)^{-2}$.
Since
$$h^2(\mu(r)(2r)^{-2})'=h^2\mu'(r)(2r)^{-2}-2^{-1}h^2r^{-3}\mu(r)>-h^2r^{-3}\mu(r),$$
we can use Lemmas 2.4 and 2.5 instead of Lemma 2.3 to conclude that the inequality (\ref{eq:5.2}) remains valid for $r\ge r_0$
with $E/2$ replaced by $E/3$.
\eproof

\end{document}